\documentclass[12pt,twocolumn]{article}

\usepackage{graphicx}
\usepackage[usenames, dvipsnames]{xcolor}
\usepackage[authoryear,round]{natbib}
\usepackage{tikz}
\usepackage{amsmath,amsfonts}
\usepackage{stmaryrd}
\usepackage[all]{xy}
\usepackage{dsfont}
\usepackage{tabularx}
\usepackage[ruled,lined,commentsnumbered]{algorithm2e}
\usepackage{graphicx}
\usepackage{booktabs}
\newcommand{\argmin}[1]{\underset{ #1}{\operatorname{argmin}}}
\newcommand{\Fro}[1]{\left\Vert #1 \right\Vert_{F}} 
\newcommand{\normeun}[1]{\left\Vert #1 \right\Vert_1} 
\newcommand{\Perm}{\mathbb{P}_p(\mathbb{R})}
\newcommand{\Tri}{\mathbb{T}_p(\mathbb{R})}
\newcommand{\pen}{\operatorname{pen}}
\newcommand{\abs}[1]{\left| #1 \right|} 

\newtheorem{djoe}{Example}
\newtheorem{hjp}{Proposition}[section]

\newtheorem{gwe}{Theorem}

\usepackage[hidelinks]{hyperref}

\title{Inferring large graphs using $\ell_1$-penalized likelihood}%
\author{Magali Champion 
  \\
Laboratoire MAP5, Universit\'e Paris Descartes,\\
Sorbonne Paris Cit\'e, France\\
    Victor Picheny \\
  INRA, UR875 Applied Mathematics and Computer Science Unit, \\
  Castanet-Tolosan, France
 \\
    and \\
    Matthieu Vignes \\
    Institute of Fundamental Sciences, Massey University,\\
     Palmerston North, New Zealand.}
    \date{}
\begin{document}

\maketitle

\begin{abstract}
We address the issue of recovering the structure of large sparse directed acyclic graphs from noisy observations of the system. 
We propose a novel procedure based on a specific formulation of the $\ell_1$-norm regularized maximum likelihood, 
which decomposes the graph estimation into two optimization sub-problems: topological structure and node order learning. 
We provide convergence inequalities for the graph estimator, as well as an algorithm to solve the induced optimization
problem, in the form of a convex program embedded in a genetic algorithm.
We apply our method to various data sets (including data from the DREAM4 challenge)
and show that it compares favorably to state-of-the-art methods.
This algorithm is available on CRAN as the \textsf{R}\/ package \texttt{GADAG}.

\end{abstract}

\noindent%

\vfill

\section{Introduction}
\label{intro}
Revealing the true structure of a complex system is paramount in many fields to identify system regulators, predict its behavior or decide where interventions are needed to disentangle direct relationships \citep{newman2003, souma2006, verma2014}. 
This problem can often be seen as a graph inference problem. Given observational data, we aim at predicting the presence (or absence) of edges between elements of the system, which form the vertices of a graph. Edges specify the relational structure of the system depicted as a graph or network. As a motivating problem, the reconstruction of Gene Regulatory Networks (GRN), which model activation and inhibition relationships between genes, is one of the main challenges in modern computational biology \citep{barabasi2004}.

A popular approach consists in assuming that the data are generated by a Directed Acyclic Graph (DAG) \citep{Pearl09}. 
DAGs are made of a collection of vertices, which stand for variables, and directed edges to model the dependency structure among the variables, avoiding self-loops and cycles.
However, inferring a DAG is a rather challenging problem. 
Firstly, the number of nodes $ p $ of the graph may
be so large that exploring relevant DAG topologies is simply infeasible, since the number of possible DAG structures is super-exponential in $ p $ \citep{Robinson73,koivisto2004, Tsamardinos06, grzegorczyk2008}. 
Another dimension flaw occurs when $ p $, even being reasonable, is larger than the number of observations, and parameter estimation is jeopardized. 
High-dimensional statistical techniques are then needed to overcome this issue \citep{Buhlmann11, giraud2014}.
Secondly, even if the ratio between $ p $ and the sample size $ n $ is not impeding model estimation, the 
nature of the data can be an additional obstacle \citep{ellis2008, guyon2010, fu2013}.
The available observational data are in general not sufficient to identify the true underlying DAG,
and can only determine an equivalence class of DAGs \citep{Pearl91}. 
This approach relies on the assumption that the joint distribution is Markov and faithful with respect to the true graph
 \citep{Spirtes00}.

A large number of methods have been proposed for estimating DAGs, including for instance score-based methods (Bayesian score, \citealt{friedman2003} or Bayesian Information Criterion, \citealt{Schwarz78}), complex space sampling \citep{zhou2011} or the PC algorithm \citep{Spirtes00}. 
The latter has been proved to be uniformly consistant in the high-dimensional case, but requires a test of conditional independences that quickly becomes computationally intractable \citep{Kalisch07}.
Recent works stress the limitations of the absence of cycles in DAGs in the study of complex systems \citep{DeSmet10}. \citet{Wright21} already described more general directed dependencies between variables when introducing genetic path analysis, but the data were then very limited. Structural Equation Modelling later introduced the notion of noise measurement \citep{Hoyle95} and \citet{Pearl09} extended them beyond linearity. Moreover, the directed cyclic graph \citep{Spirtes95} framework received little attention as compared to its acyclic counterpart. Finally, the actual discovery of causal cycles requires temporal data, \textit{e.g.} in the context of dynamic Bayesian networks \citep{Perrin03,Dondelinger13}, data which are difficult and very expensive to collect despite efforts in this direction \citep{Sachs09}.
 
In this work, we focus on Gaussian structural equation models
associated with maximum likelihood estimators (MLE).
In the last years, the $\ell_0$-regularization of the MLE drew the attention of a large number of works since it leads to infer sparse graphs. In DAGs, a not necessarily topological ordering of the nodes can always be defined according to edge distribution \citep{kahn1962}. Identifying this ordering is known to be a challenging problem \citep{Cook85}. 
Additional data, like gene knock-out data or more general perturbations data \citep{Maathuis10,Shojaie14} can give information in that way.
More generally, biological prior knowledge, retrieved from specific data bases, can be used to assist the network reconstruction algorithm \citep{Husmeier07}, or a partial knowledge of the network can inform the inference process efficiently, \textit{e.g.} in the semi-supervised framework of \citet{Mordelet08}.
For a known order among the variables in the graph, \citet{Shojaie10} present results for the estimation of high-dimensional graphs based on independent linear regressions using an adaptive Lasso scheme.
When the order of the variables is unknown, \citet{vandeGeer13} studied the convergence of the $\ell_0$-penalized likelihood. However, the $\ell_0$-regularized approaches \citep{Silander06,Hauser12} remain impractical for estimating graphs with more than 20 vertices, either due to an exhaustive exploration of the set of DAGs or overfitting \citep{Chen08}. Quite recently, \citet{Aragam15} explored a penalized least-squares estimator for a variety of concave penalization terms. They obtain theoretical guarantees of sparsity bounds and consistency results in model selection. Their theoretical results would greatly benefit an implementation of the methods and an empirical study to demonstrate the effectiveness of the approach. The unifying framework for pseudolikelihood-based graphical modelling of \citet{Khare15} extends classical regularization methods. The authors obtain theoretical convergence results, and offer an algorithm with an associated implementation. Their simulation results are quite promising, in particular in terms of computational time.

Our objective consists in overcoming this drastic dimensional limitation, and find inference strategies for graphs with up to several hundred nodes.
Such strategies must ensure a high level of sparsity, be supported by computationally affordable algorithms, 
while preserving sound theoretical bases.
Here, we propose to use the $\ell_1$-regularization, similarly to \citet{fu2013} and \citet{Shojaie10}, to penalize the MLE. From a computational point of view, this regularization makes the criterion to maximize partially convex 
while ensuring sparse estimates. Our contribution is two-fold: firstly, we provide convergence inequalities that guarantee good theoretical performances of our proposed estimator in the sparse high-dimensional setting. Secondly, we provide an efficient algorithm to infer the true unknown DAG, in the form of a convex program embedded in a genetic algorithm.

The next section covers the model definition and the associated penalized MLE problem. Section \ref{sec:convergence} details the convergence inequalities, and Section \ref{sec:algorithm} our inference algorithm. 
Section \ref{sec:experiments} reports numerical experiments both on toy problems and realistic data sets.

\section{The $\ell_1$-penalized likelihood for estimating DAGs}\label{sec:likelihood}
\subsection{DAG's modelling and estimation}\label{sec:model}
This work considers the framework of an unknown DAG $\mathcal{G}_0=(V,E)$, consisting of vertices $V=\{1,...,p\}$ and a set of edges $E\subseteq V\times V$.
The $p$ nodes are associated to random variables $X^1,...,X^p$.
A natural approach, developped by \citet{Meinshausen06} to solve the network inference problem is to consider that each variable $X^i$ ($1\leq i\leq p$) of the DAG can be represented as a linear function of all other variables $X^j$ ($j\neq i$) through the Gaussian Structural Equation Model:
\begin{equation}\label{eq:model}
\forall j \in \llbracket 1, p \rrbracket, \ \ X^j = \sum_{i=1}^p (G_0)_i^j X^i + \varepsilon^j,
\end{equation}
with $\varepsilon^j \sim \mathcal{N}(0,\sigma^2_j)$ 
($\sigma^2_j$ known)
a Gaussian residual error term. The set of edges $E$, which is assumed to be of size $s$ ($s\leq p(p-1)/2$), corresponds to the  non-zero coefficients of $G_0$, \textit{i.e.} $(G_0)^j_i$ encodes the relationship from variable $X^i$ to variable $X^j$. 

Assume that we observe an $n$-sample consisting of $n$ i.i.d. realizations $(X^1,...,X^p)$ of Equation (\ref{eq:model}), distributed according to a $\mathcal{N}(0,\Sigma)$ law where $\Sigma$ is non-singular. We denote by $X:=(X^1,...,X^p)$ the $n\times p$ data matrix. The relations between the variables can then be represented in its matrix form:
\begin{equation}\label{eq:modelmatrix}
X=XG_0 + \varepsilon,
\end{equation}
where $G_0=((G_0)_i^j)_{1\leq i,j\leq p}$ is the $p\times p$ matrix compatible with the graph $\mathcal{G}_0$ and $\varepsilon:=(\varepsilon^1,...,\varepsilon^p)$ is the $n\times p$ matrix of noise vectors.

The negative log-likelihood of the model is then \citep{Rau13}:
\begin{multline}\label{prop:lv}
\ell(G)= \frac{np}{2} \log (2\pi) +
 n\sum_{j=1}^p \log \sigma_j \\+ \sum_{k=1}^n \sum_{j=1}^p \frac{1}{\sigma_j^2} \left( X_k(I-G)^j \right)^2,
\end{multline}
where $I$ denotes the $p\times p$ identity matrix and $X_k$ the vector of length $p$ corresponding to the $k$-th observation of $X^1,...,X^p$.

To recover the structure of the DAG $\mathcal{G}_0$ and make the estimated graph sparse enough, we focus on a penalized maximum likelihood procedure \citep{bickel2006}:
\begin{equation}\label{eq:lv}
\hat{G} = \argmin{G\in \mathcal{G}_{\mbox{\scriptsize{DAG}}}}\ \  \{ \ell(G) + \lambda \pen(G) \},
\end{equation}
where $\ell(.)$ is the negative log-likelihood of Equation (\ref{prop:lv}), $\pen(.)$ is a penalization function, 
$\lambda$ is a trade-off parameter between penalization and fit to the data, and $\mathcal{G}_{\mbox{\scriptsize{DAG}}}$ is the set of $ p \times p $ matrices compatible with a DAG over $ p $ nodes.

In the setting of
Gaussian structural equation models with equal noise variance, \citet{Peters11, peters2014} showed that the true DAG was identifiable for respectively discrete and continuous data. In a nutshell, it implies that the true DAG could be inferred, not just the Markov equivalence class of the underlying DAG - a partially directed graph exactly encoding the conditional dependency structure.
Using an $\ell_0$-norm regularization in Equation (\ref{eq:lv}) is an attractive option to infer sparse graphs. 
From a computational point of view, the main difficulty when solving the optimization problem in Equation (\ref{eq:lv}) lies in exploring the set of DAGs $\mathcal{G}_{\mbox{\scriptsize{DAG}}}$. \citep{Chickering96} showed it to be an
NP-hard problem:
an $\ell_0$-regularization does not set a favorable framework for this task.
To avoid the whole exploration of $\mathcal{G}_{\mbox{\scriptsize{DAG}}}$, a dynamic programming method has been proposed in \citet{Silander06}, using a particular decomposition of the $\ell_0$-penalized maximum likelihood. 
The greedy equivalent search algorithm of \citet{Chickering02} is a hill climbing alternative method. \citep{Hauser12} rather restricted the search space to the smaller space of equivalence classes, and they provide 
an efficient algorithm without enumerating all the equivalent DAGs. They showed that they are asymptotically optimal under a faithfulness assumption (\textit{i.e.} independences in the distribution are those read from $ \mathcal{G}_{0} $).
However, none of the approaches above can be used on high-dimensional data to estimate graphs with a large number of nodes.
In this context, we focus on the $\ell_1$-norm convex regularization instead of $\ell_0$
for its sparse, high-dimensional and computational properties. 

The $\ell_1$-regularization clearly improves the computation of (\ref{eq:lv}). It allows us to write a convex formulation of the problem
 (see Section \ref{sec-ref}).
Given Equation (\ref{prop:lv}) and omitting constant terms, the $\ell_1$-penalized likelihood estimator we consider is:
\begin{equation}\label{eq:lv2}
\hat{G} = \argmin{G\in \mathcal{G}_{DAG}} \  \ \left\{ \frac{1}{n} \Fro{X(I-G)}^2 + \lambda \normeun{G} \right\},
\end{equation}
where, for any matrix $M:=(M_i^j)_{1\leq i,j \leq p}$, 
we denote by $\Fro{M} = \sum_{i,j} ( M_i^j )^2$ the Frobenius norm and by $\normeun{M} = \sum_{i,j} | M_i^j|$ the $\ell_1$-norm.

\subsection{A new formulation for the estimator}\label{sec-ref}
We propose here a new formulation of the minimization problem of Equation (\ref{eq:lv2}). It 
naturally decouples the initial problem into two steps of the minimisation procedure: node ordering and graph topology search.
A key property is that any DAG leads to a topological ordering of its vertices, denoted $ \leq $, where a directed path from node $ X^i $ to node $ X^j $ is equivalent to $X^j \leq X^i$ \citep{kahn1962, Cormen01} (see Example \ref{ex:DAG} below for more explanations). This ordering is not unique in general.
Proposition \ref{prop:decomp} from \citet{Buhlmann13} then gives an equivalent condition for a matrix to be compatible with a DAG.
\begin{hjp}[B\"{u}hlmann, 2013]\label{prop:decomp}
A matrix $G$ is compatible with a DAG $\mathcal{G}$ if and only if there exists a permutation matrix $P$ and a strictly lower triangular matrix $T$ such that:
$$G=PTP^T.$$
\end{hjp}
Graphically, the permutation matrix sets an ordering of the nodes of the graph and is associated to a complete graph. The strictly lower triangular matrix $T$ sets the graph structure, \textit{i.e.} the non-zero entries of $G$, as illustrated in Example \ref{ex:DAG}.

\begin{djoe}\label{ex:DAG}
Consider the DAG $\mathcal{G}$ given in Figure \ref{fig:DAG} (left). The corresponding matrix $G$ can then be written as the strictly lower-triangular matrix $T$ by permutation of its rows and columns using $P$:
$$
G = \begin{pmatrix}
0&0&0&7&5\\
2&0&1&6&4\\
0&0&0&0&0\\
0&0&0&0&3\\
0&0&0&0&0
\end{pmatrix}=PTP^T,$$

$$
T=\begin{pmatrix}
0&0&0&0&0\\
0&0&0&0&0\\
3&0&0&0&0\\
5&0&7&0&0\\
4&1&6&2&0
\end{pmatrix}\ \ \mbox{and} \ \ 
P=\begin{pmatrix}
0&0&0&1&0\\
0&0&0&0&1\\
0&1&0&0&0\\
0&0&1&0&0\\
1&0&0&0&0
\end{pmatrix}.$$
Looking at the non-zero values of $P$ column by column, $P$ defines a node hierarchy $X^5 \leq X^3\leq X^4 \leq X^1 \leq X^2$ compatible with the topological orderings of $\mathcal{G}$.
Graphically, $P$ is associated to the complete graph represented in Figure \ref{fig:DAG} (bottom). The dashed edges then correspond to the lower zero entries of $T$.
Note that since $X^3$ is not connected with $X^1$, $X^4$ and $X^5$, four topological ordering are possible ($X^5\leq X^4 \leq X^1 \leq X^3$, $X^5\leq X^4 \leq X^3 \leq X^1$, $X^5 \leq X^3 \leq X^4 \leq X^1$ and $X^3 \leq X^5\leq X^4 \leq X^1$).

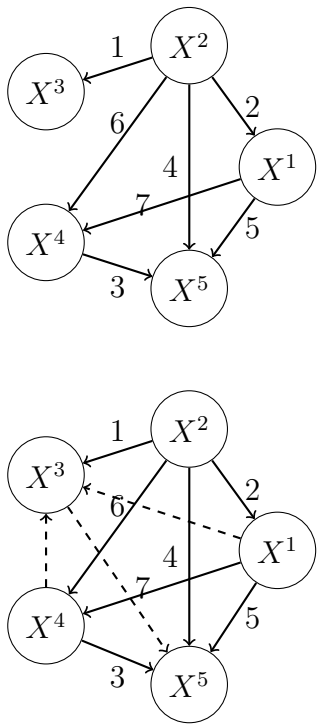
\begin{figure}[!ht]
\begin{center}
\begin{tikzpicture}[scale=1.7]
\node[circle, draw=black] (X1bis) at (4,0) {$X^1$};
\node[circle, draw=black] (X2bis) at (3.31,0.95) {$X^2$};
\node[circle, draw=black] (X3bis) at (2.19,0.59) {$X^3$};
\node[circle, draw=black] (X4bis) at (2.19,-0.59) {$X^4$};
\node[circle, draw=black] (X5bis) at (3.31,-0.95) {$X^5$};
\draw[->,thick] (X1bis)--(X4bis) node[midway,above,left]{7};
\draw[->,thick] (X2bis)--(X3bis) node[midway,above]{1};
\draw[->,thick] (X2bis)--(X5bis) node[midway,left]{4};
\draw[->,thick] (X2bis)--(X4bis) node[midway,left,above]{6};
\draw[->,thick] (X1bis)--(X5bis) node[midway,right]{5};
\draw[->,thick] (X4bis)--(X5bis) node[midway,below]{3};
\draw[->,thick] (X2bis)--(X1bis) node[midway,right]{2};

\node[circle, draw=black] (X1) at (4,-3) {$X^1$};
\node[circle, draw=black] (X2) at (3.31,-2.05) {$X^2$};
\node[circle, draw=black] (X3) at (2.19,-2.41) {$X^3$};
\node[circle, draw=black] (X4) at (2.19,-3.59) {$X^4$};
\node[circle, draw=black] (X5) at (3.31,-4.05) {$X^5$};
\draw[->,thick] (X2)--(X5) node[midway,left]{4};
\draw[->,thick] (X2)--(X4) node[midway,left,above]{6};
\draw[->,thick] (X2)--(X1) node[midway,right]{2};
\draw[->,thick] (X1)--(X4) node[midway,above,left]{7};
\draw[->,thick] (X2)--(X3) node[midway,above]{1};
\draw[->,thick] (X1)--(X5) node[midway,right]{5};
\draw[->,thick] (X4)--(X5) node[midway,below]{3};
\draw[->,thick,dashed] (X4)--(X3);
\draw[->,thick,dashed] (X1)--(X3);
\draw[->,thick,dashed] (X3)--(X5);
\end{tikzpicture}
\end{center}
\caption{An example of DAG $\mathcal{G}$ (top) and the action of $P$ and $T$ on $\mathcal{G}$: $P$ is associated to a complete graph that orders the nodes of the graph (bottom) and $T$ sets the weights on the edges. The dashed edges correspond to null weight edges (a zero entry in $T$).}\label{fig:DAG}
\end{figure}
\end{djoe}

Using Proposition \ref{prop:decomp}, the estimator in (\ref{eq:lv2}) leads to the following equivalent optimization problem:
\begin{multline}\label{eq:PT}
(\hat{P},\hat{T}) =\\
 \argmin{(P,T) \in \mathcal{C}} \ \ \left\{ \frac{1}{n} \Fro{X(I-PTP^T)}^2 + \lambda \normeun{T} \right\},
\end{multline}
where the optimization space of constraints $\mathcal{C}$ is defined as $\mathcal{C} = \Perm \times \Tri$, with $\Perm$ the set of permutation matrices and $\Tri$ the set of strictly lower-triangular matrices.
Note that a similar formulation has already been proposed by \citet{vandeGeer13} to ensure good theoretical properties for the $\ell_0$-penalized log-likelihood estimation. However, it has never been exploited from a computational point of view to recover the graph structure optimizing problem (\ref{eq:lv2}).
In the following two sections, we propose a theoretical analysis of the proposed estimator (Section \ref{sec:convergence}) and a computationally effficient algorithm to solve Problem (\ref{eq:PT}) (Section \ref{sec:algorithm}).

\section{Convergence inequalities for the DAG estimation}\label{sec:convergence}
The main result of this section deals with convergence rates: in Theorem \ref{th:main}, we provide upper bound for error associated with the $\ell_1$-penalized maximum likelihood estimator considered in Equation (\ref{eq:PT}), both in prediction (Equation \ref{eq:pred}) and estimation (Equation \ref{eq:inegest}). 
Following the works of \citet{vandeGeer13} on the $\ell_0$-penalized maximum likelihood estimator and of \citet{Bickel09} on the Lasso and the Dantzig Selector, we obtain two convergence results under some mild sparsity assumptions, when the number of variables is large but upper bounded by a function $\varphi(n) $ of the sample size $ n $.

\subsection{Estimating the true order of variables}
For a known ordering among the variables of the graph \citep{Shojaie10}, an unrealistic assumption in many applications, the DAG inference problem can be cast in a convex optimization problem.
To provide convergence inequalities of the proposed estimator in the most general case of an unknown order we consider here, we first focus on the problem of estimating the true variable order.
Let us denote by $\Pi_0$ the set of permutation matrices compatible with the true DAG $\mathcal{G}_0$:
$$\Pi_0=\left\{ P\in \Perm, P^TG_0P \in \Tri \right\}.$$
$\Pi_0$ contains one or more permutation matrice(s) (see Example \ref{ex:DAG}). We will have to make a decision as to whether the estimated order of variables $\hat{P}$ given by Equation (\ref{eq:PT}) is in $\Pi_0$ or not. 

To answer this question, we investigate the effect of learning an erroneous order of variables $P\notin \Pi_0$. We introduce the following additional notations:
for any permutation matrix $P\in \Perm$, we denote by $G_0(P)$ the matrix defined as:
$$ G_0(P)=PT_0P^T,$$
with $T_0 = P_0^T G_0 P_0$ a lower triangular decomposition of $G_0$. From a graphical point of view, while $P\notin \Pi_0$, the graph $\mathcal{G}_0(P)$ associated to $G_0(P)$ is obtained from $\mathcal{G}_0$ by permuting some of its nodes (see Example \ref{ex:G0(P)}), otherwise, if $P\in \Pi_0$, $\mathcal{G}_0(P)=\mathcal{G}_0$.
We also denote by $\varepsilon(P):=X-XG_0(P)$ the associated residual term.
We denote by $\Omega(P)$ the covariance matrix of $\varepsilon(P)$ and $\omega_j(P):=\operatorname{Var}(\varepsilon^j(P))$ the associated noise variances of each node.

With these notations and checking that the assumptions presented in Section \ref{sec-hyp} hold, we ensure that, with large probability, we choose a right order of variables and the estimated graph converges to the true graph when $ n $ and $ p $ grow to infinity (see Section \ref{sec-main}).

\begin{djoe}\label{ex:G0(P)}
Let $$P=\begin{pmatrix}
0&0&0&0&1\\
0&0&0&1&0\\
0&1&0&0&0\\
1&0&0&0&0\\
0&0&1&0&0
\end{pmatrix}
\notin \Pi_0$$ a wrong permutation.

In Figure \ref{fig:permgraph}, we represent the permuted graph $\mathcal{G}_0(P)$ (bottom) associated to the graph $\mathcal{G}_0$ (top). 
The latter is obtained from $\mathcal{G}_0$ after permutation of its nodes using $PP_0^T$, where $P_0$ (corresponding to the matrix $P$ in Example \ref{ex:DAG}) defines a right order of variables.
\begin{figure}[!ht]
\begin{center}
\begin{tikzpicture}[scale=1.7]
\node[circle, draw=black] (X1) at (4,0) {$X^1$};
\node[circle, draw=black] (X2) at (3.31,0.95) {$X^2$};
\node[circle, draw=black] (X3) at (2.19,0.59) {$X^3$};
\node[circle, draw=black] (X4) at (2.19,-0.59) {$X^4$};
\node[circle, draw=black] (X5) at (3.31,-0.95) {$X^5$};
\draw[->,thick] (X3)--(X4) node[midway,left]{7};
\draw[->,thick] (X2)--(X3) node[midway,above]{1};
\draw[->,thick] (X2)--(X5) node[midway,left]{4};
\draw[->,thick] (X2)--(X4) node[midway,left,above]{6};
\draw[->,thick] (X3)--(X5) node[midway,left,below]{5};
\draw[->,thick] (X5)--(X1) node[midway,below,right]{3};
\draw[->,thick] (X2)--(X1) node[midway,right,above]{2};

\node[circle, draw=black] (X1bis) at (4,-3) {$X^4$};
\node[circle, draw=black] (X2bis) at (3.31,-2.05) {$X^1$};
\node[circle, draw=black] (X3bis) at (2.19,-2.41) {$X^2$};
\node[circle, draw=black] (X4bis) at (2.19,-3.59) {$X^5$};
\node[circle, draw=black] (X5bis) at (3.31,-4.05) {$X^3$};
\draw[->,thick] (X3bis)--(X4bis) node[midway,left]{$7$};
\draw[->,thick] (X2bis)--(X3bis) node[midway,above]{$1$};
\draw[->,thick] (X2bis)--(X5bis) node[midway,left]{$4$};
\draw[->,thick] (X2bis)--(X4bis) node[midway,left,above]{$6$};
\draw[->,thick] (X3bis)--(X5bis) node[midway,left,below]{$5$};
\draw[->,thick] (X5bis)--(X1bis) node[midway,below,right]{$3$};
\draw[->,thick] (X2bis)--(X1bis) node[midway,right,above]{$2$};

\end{tikzpicture}
\end{center}
\caption{The graph $\mathcal{G}_0$ (top) and the permuted graph $\mathcal{G}_0(P)$ (bottom) associated to the permutation $P$.}\label{fig:permgraph}
\end{figure}
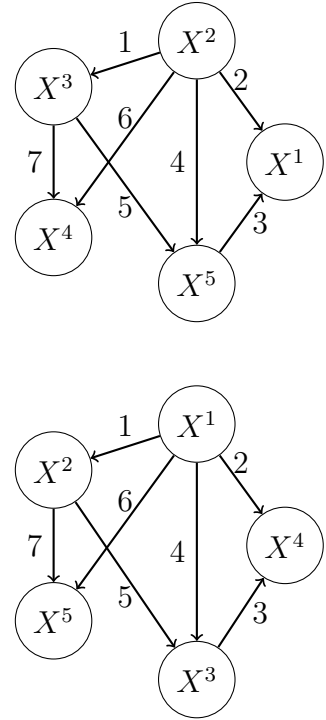

\end{djoe}

\subsection{Assumptions on the model}\label{sec-hyp}
For a square matrix $M\in \mathcal{M}_{p\times p}(\mathbb{R})$ and a subset $\mathcal{S}$ of $\llbracket 1,p \rrbracket^2$, 
we denote by $M_{\mathcal{S}} \in \mathcal{M}_{p\times p}(\mathbb{R})$ the matrix that has the same elements as $M$ on $\mathcal{S}$ and zero on the complementary set $\mathcal{S}^C$ of $\mathcal{S}$.
We now introduce the assumptions we used to obtain statistical properties of our estimator.

\paragraph{Hypotheses} \textcolor{white}{Nothing}

\noindent $\mathbf{H_{1}}$ There exists $\sigma^2>0$ such that $$\forall j \in \llbracket 1,p \rrbracket, \operatorname{Var}(\varepsilon^j) = \sigma^2.$$

\noindent $\mathbf{H_{2}}$ There exists $\sigma_0^2$, independent of $p$ and $n$, such that
$$\max_{1 \leq j\leq p} \operatorname{Var}(X^j) \leq \sigma_0^2.$$

\noindent $\mathbf{H_{3}}$ There exists $\lambda^*>0$ such that the minimal eigenvalue of the covariance matrix $\Sigma$ of $X$ satisfies $$\lambda_{min} \geq \lambda^* >0.$$

\noindent $\mathbf{H_{4}}$ There exists $g_{max} < \infty$ such that the maximal weight of the DAG $\mathcal{G}_0$ is bounded  $$\displaystyle \max_{1\leq i,j\leq p} \abs{(G_0)_i^j} \leq g_{max}.$$ 

\noindent $\mathbf{H_{5}}$ The number of nodes $p$ satisfies
$$p \log p = \mathcal{O}(n).$$

\noindent $\mathbf{H_{6}}$ There exists $\kappa(t)>0$ with $1\leq t \leq p^2$ such that:
$$\min \left\{ \frac{\Fro{XM}}{\sqrt{n} \Fro{M_{\mathcal{S}}}} \right\}
\geq \kappa(t),$$
where the minimum is taken over the set of $p\times p$ matrices satisfying $\normeun{M_{\mathcal{S}^C}}\leq 3 \normeun{M_{\mathcal{S}}}$, with $\mathcal{S} \subset \llbracket 1,p \rrbracket^2$ and $\abs{\mathcal{S}} \leq t$.

\noindent $\mathbf{H_{7}}$
There exists 
$0<\eta \leq C\frac{n}{p\log p}\times\frac{1}{ \sqrt{s}}$
 such that, for all permutations $P\notin \Pi_0$,
$$\frac{1}{p} \sum_{j=1}^p \left( \abs{\omega_j(P)}^2 -1\right)^2 > \frac{1}{\eta}.$$
%

Assumption $\mathbf{H_1}$ states that the noise variances are the same among all variables. This assumption is clearly hard to test in practice but makes the problem identifiable and ensures that we can recover the true DAG. Otherwise, minimizing (\ref{eq:lv2}) only leads to the identification of one element of the Markov equivalence class of the true DAG (partially directed graph). To simplify the theoretical results and proofs, until the end of this work, we assume that the noise variances $\sigma^2$ are equal to one. Our results are still valid even if $\sigma^2\neq 1$, by small modifications in the constant terms as long as they are all equal.

 Assumption $\mathbf{H_{5}}$ deserves a special attention since it bounds the high dimensional setting. The considered problem is obviously non-trivial and requires a sufficient amount of information. A more detailled discussion about assumptions $\mathbf{H_3}$ and $\mathbf{H_5}$ is proposed in Section \ref{sec:cond}.

Assumption $\mathbf{H_{6}}$ is a natural extension of the Restricted Eigenvalue condition of \citet{Bickel09} to our multi-task setting. More precisely, denoting 
\begin{center}
\begin{tikzpicture}[scale=1]
\draw (0,0) node{$\tilde{X}=\left(
     \raisebox{0.5\depth}{%
       \xymatrixcolsep{1ex}%
       \xymatrixrowsep{1ex}%
       \xymatrix{
         X \ar @{-}[rrdd] & & 0\\
         & & \\
         0 & & X
       }%
     }
   \right)$};
   \draw[<->] (-0.4,-1)--(1.3,-1);
   \draw[<->] (1.7,-0.8)--(1.7,0.8);
   \draw (0.45,-1) node[below]{$p^2$};
   \draw (1.7,-0) node[right]{$n\times p$,};
\end{tikzpicture}
\end{center}
$\mathbf{H_{6}}$ is equivalent to assuming that the Gram matrix $\frac{\tilde{X}\tilde{X}^T}{n}$ is non-degenerate on a restricted cone \citep{Lounici09,Buhlmann11}. Notice that this condition is very classical in the literature. It yields good practical performance even for small sample sizes, and some recent works discuss an
accurate population eigenvalue estimation
even in a large dimension setting \citep{mestre2008, elkaroui2008, liu2014, ledoit2015}.

The last assumption $\mathbf{H_{7}}$ is an identifiability condition needed to ensure that the estimated permutation $\hat{P}$ is in $\Pi_0$. This assumption was introduced by \citet{vandeGeer13} as the ``omega-min'' condition. In a sense, it separates the set of compatible permutations from its complement in a finite sample scenario.

\subsection{Main result}\label{sec-main}
The result we establish in this section is double-edged: (a) with large probability, we ensure that the estimated $ \hat{P} $ belongs to $\Pi_0$, and (b) we provide convergence inequalities both in prediction and estimation for the graph estimated from the minimisation problem (\ref{eq:PT}). This result clearly states the desirable theoretical properties of the derived estimator, assuming reasonable conditions on the complex system embedding the data.

\begin{gwe}\label{th:main}
Assume that $\mathbf{H_{1-7}}$ are satisfied, with $s$ $\subset \llbracket 1,p^2 \rrbracket$ in $\mathbf{H_6}$ such that $\sum_{i,j} \mathds{1}_{(G_0)_i^j \neq 0} \leq s$ ($G_0$ is $s$-sparse). Let $\lambda = 2C\sqrt{s^{1/2}\frac{\log p}{n}}.$ Then, with probability greater than $1-5/p$, any solution $\hat{G}=\hat{P} \hat{T} \hat{P}^T$ of the minimization problem (\ref{eq:PT}) satisfies that $\hat{P}\in \Pi_0$. Moreover, with at least the same probability, the following inequalities hold:
\begin{equation}\label{eq:pred}
\frac{1}{n} \Fro{X\hat{G}-XG_0}^2 \leq  \frac{16C^2}{\kappa^2(s)} s^{3/2} \frac{\log p }{n}.
 \end{equation}

 \begin{equation}\label{eq:inegest}
\normeun{\hat{G}-G_0} \leq \frac{16 C}{\kappa^2(s)} \sqrt{s^{5/2} \frac{\log p}{n}}. 
 \end{equation}
\end{gwe}

The proof of this result is deferred in Section C of the Supplementary Materials.

Theorem \ref{th:main} states that with probability at least $1-5/p$, we choose a compatible order of variables over the set of permutations. 
Inequalities (\ref{eq:pred}) and (\ref{eq:inegest}) give non-asymptotic upper bounds on the loss under conditions depending on $s$, $p$ and $ n $ (see Section \ref{sec:cond}). They also ensure that the estimated $\hat{T}$ is close to the true $T_0$ with large probability.

\subsection{Discussion on the high-dimensional scenario} \label{sec:cond}
\paragraph{Sparsity of the graph}
Assumption $\mathbf{H_7}$ and Theorem \ref{th:main} naturally require a trade-off between signal sparsity, dimensionality and sample size. In the ultra sparse 
regime (where the sparsity $s$ of the true graph is bounded by $s^*>0$), Theorem \ref{th:main} provides convergence inequalities for $\hat{G}$ choosing $\eta\leq \frac{\alpha}{\sqrt{s^*}}$ with $p\log(p) = \alpha n$ in Assumption $\mathbf{H}_7$.

In the standard sparsity scenario, if $s$ is at least of the order of $p$, then $\eta$ should be of the order of $\alpha / \sqrt{p}$, which is unrealistic as $p\rightarrow + \infty$. This case thus requires a stronger dimensional assumption $\mathbf{H}_5$. Taking at least $p^{2} \log(p) = \mathcal{O}(n)$ ensures a good estimation of the graph. 

Note however that universal conditions cannot be overcome and the ultra-high dimension settings (\textit{e.g.} \citet{wainwright2009, verzelen2012}) is an insurmountable limit.

\paragraph{Minimal eigenvalue condition}
Assumption $\mathbf{H_3}$ ensures that the minimal eigenvalue of the covariance matrix $\Sigma$ of $X$ is not too small. In the high-dimensional scenario, this could be hard to verify, $\lambda_{min}$ decreasing while $n,p$ growing to infinity \citep{Hogben07}. A natural bound for $\lambda_{min}$ is:
\begin{equation}\label{eq:lam}
\lambda_{min} \geq \frac{1}{p\max\left(1,g_{max}^2\right) (1+\sqrt{s})},
\end{equation}
with $g_{max}$ and $s$ as in $\mathbf{H_{4}}$ and Theorem \ref{th:main}.

Assumption $\mathbf{H_3}$ can thus be relaxed by allowing $\lambda_{min}$ to decrease with $1/p\sqrt{s}$. The price to pay for this relaxation is a data dimensionality reduction $p^3 \log (p) = \mathcal{O}(n)$, which automatically implies:
$$\frac{3\lambda_{min}}{4} -2 \sqrt{\frac{\log (p)}{n}} -3\sigma_0 \sqrt{\frac{2p\log (p)}{n}} >0,$$ with Equation (\ref{eq:lam}) (for more details, see Section A of the Supplementary Material \textit{Proof details}). 

\section{Inference algorithm}\label{sec:algorithm}

\subsection{Global algorithm overview}
In this section, we present GADAG (\textit{Genetic Algorithm for learning Directed Acyclic Graphs}), a computational procedure devoted to solve Equation (\ref{eq:PT}) and available as a \texttt{R} package on CRAN at \url{https://cran.r-project.org/package=GADAG}. 
Although decomposing the original problem made it more natural to handle, this problem is still a very challenging task
from an optimization point of view, due to the different nature of the variables $P$ and $T$, the non-convexity of the cost function
and the high dimension of the search space.

An intuitive approach consists in using an alternating scheme: 
one of the variables $ P $ or $ T $ is fixed and the other one is sought so as to optimize the score function,
then the roles of $ P $ and $ T $ are reversed and the procedure is repeated iteratively until convergence for some criterion \citep{Csiszar84}.
However, the structure of our problem does not allow us to use such a scheme: looking for an optimal $T$ given a fixed $P$ makes sense, 
but changing $P$ for a fixed $T$ does not. 

\begin{figure}[!ht]
\begin{center}
 \begin{tikzpicture}[scale=0.55]
 \draw (-4.5,8.5) rectangle (4.5,9.5);
  \draw (-4.5,6.5) rectangle (4.5,7.5);
  \draw (-4.5,4.5) rectangle (4.5,5.5);
\draw[>=stealth,->](0,10.5)--(0,9.5);
  \draw[>=stealth,->](0,8.5)--(0,7.5);
  \draw[>=stealth,->](0,6.5)--(0,5.5);
   \draw[>=stealth,->](0,4.5)--(0,3.5);

   \draw[-](8,10.5)--(0,10.5);
    \draw[-](8,10.5)--(8,2.5);
     \draw[-](8,2.5)--(4,2.5);
     \draw[>=stealth,->](0,1.5)--(0,0.5);

   \draw(0,9) node{Choice of $P$};
 \draw(0,7) node{Search of an optimal $T^*$};
\draw(0,5) node{Evaluate the likelihood};  
 \draw(0,2.5) node{Problem solved?};
   \draw[-](0,3.5)--(4,2.5);
   \draw[-](0,3.5)--(-4,2.5);
   \draw[-](4,2.5)--(0,1.5);
   \draw[-](-4,2.5)--(0,1.5);

      \draw(0,0) node{END};
      \draw(1,1) node{YES};
    \draw(6,2) node{NO};   
 \end{tikzpicture}
 \end{center}
 \caption{Overview of our hybrid algorithm GADAG.}\label{fig:recap}
 \end{figure}
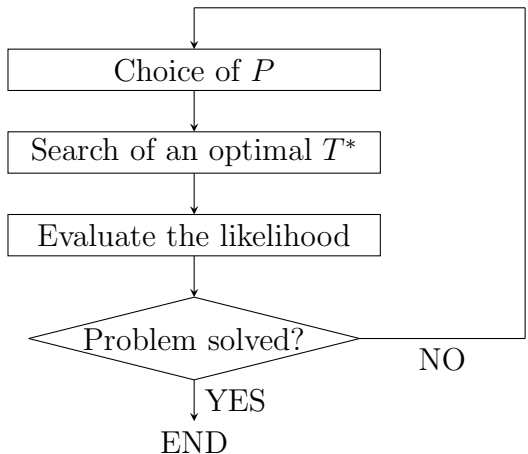

In our inference algorithm GADAG, an outer loop is used to perform the global search among the DAGs space, which is driven by the choice of $P$, while a nested loop is used to find an optimal $T$ for each given fixed $P$ (see Figure \ref{fig:recap}).
As we show in the following, population-based metaheuristics algorithms are a natural and efficient choice for exploring the space of permutation matrices (Section \ref{sec:Poptim}). 
The nested optimization problem can be resolved using a steepest descent approach (Section \ref{sec:Toptim}).

\subsection{Graph structure learning when the variable order is fixed}\label{sec:Toptim}

Assume first that the variable ordering $P\in \Perm$ is fixed. The problem of inferring a graph is then reduced to estimating the graph structure, which can be solved by finding a solution of:
\begin{equation}\label{eq:minT}
\min_{T\in \Tri} \  \ \left\{ \frac{1}{n} \Fro{X(I-PTP^T)}^2 + \lambda \normeun{T} \right\}.
\end{equation}

The minimization problem given by Equation (\ref{eq:minT}) looks like a well-studied problem in machine learning, as it is closely related to the $\ell_1$-constrained quadratic program, known as the Lasso in the statistics literature \citep{Tibshirani96}.
Indeed, the $\ell_1$-regularization leads to variable selection and convex constraints that make the optimization problem solvable. 
We note here that this allows us to always provide a locally optimal solution, \textit{i.e.} optimal weight estimates given a hierarchy between the nodes.

A large number of efficient algorithms are available for computing the entire path of solutions as $\lambda$ is varied,  \textit{e.g.} the LARS algorithm of \citet{Efron04} and its alternatives.
For example, in the context of the estimation of sparse undirected graphical models, \citet{Meinshausen06} fit a Lasso model to each variable, using the others as predictors, and define some rules for model symmetrization as they do not work on DAGs. The graphical Lasso (or glasso, \citealt{friedman2007}) algorithm directly relies on the estimation of the inverse of a structure covariance matrix assumed to be sparse. Improvements were proposed for example by \citet{duchi2008} (improved stopping criterion) and \citet{witten2011} (estimation of a block-diagonal matrix). Other authors propose to solve the optimization problem using an adaptation of classical optimization methods, such as interior point \citep{Yuan07} or block coordinate descent methods \citep{Banerjee08, friedman2007}.

We propose here an original convex optimization algorithm to find the solution in Equation (\ref{eq:minT}) in a form similar to a steepest descent algorithm. Our proposed algorithm is much quicker than a glasso approach, a desirable feature as it will run at each iteration of the global algorithm (see the ``Search of an optimal $ T^{*} $'' box in Figure~\ref{fig:recap} and the ``Evaluate the new individuals'' item in Algorithm~\ref{algo:global}). Moreover, its mechanistic components (see Section B of the Supplementary Material \textit{Proof details}) allowed us to derive the theoretical results of Theorem~\ref{th:main}.
The proposed scheme can be seen as an adaptation of the LARS algorithm with matrix arguments.

Let $(T_k)_{k\geq 0}$ the sequence of matrices defined for all $i,j \in \llbracket 1,p \rrbracket^2$ as:
\begin{equation}\label{eq:minTfin}
(T_{k+1})_i^j = \operatorname{sign} \left(  (U_k)_i^j \right) \max \left( 0, \abs{(U_k)_i^j} - \frac{\lambda}{L} \right),
 \end{equation}
\noindent where for all $ k\geq  0, U_k = T_k - \frac{ \nabla \left( \frac{1}{n} \Fro{X(I-PT_kP^T)}^2 \right)}{L}$, $L$ is the Lipschitz constant of the gradient function $\nabla \left( \frac{1}{n} \Fro{X(I-PT_kP^T)}^2 \right)$ and $\operatorname{sign}()$ is the sign of any element.
Then, a solution of (\ref{eq:minT}) is given by performing Algorithm \ref{algo:minT}, where:
\begin{itemize}
\item the projection $\operatorname{Proj}_{\Tri}(T)$ of any $ p \times p $ real-valued matrix $T=((T_k)_i^j)_{i,j}$ on the set $\Tri$ is given by
\begin{equation}\label{eq:projT}
\left( \operatorname{Proj}_{\Tri}(T_k) \right) _i^j = \begin{cases}
 0& \mbox{if $i<j$},\\
 \left( T_k \right)_i^j & \mbox{otherwise}.
 \end{cases}
 \end{equation}
 \item the gradient of $\frac{1}{n} \Fro{X(I-PT_kP^T)}^2$ is 
 \begin{multline}\label{eq:grad}
\nabla \left( \frac{1}{n} \Fro{X(I-PT_kP^T)}^2 \right)=\\
 -\frac{2}{n}  (XP)^T (X-XPT_kP^T)P.
\end{multline} 
\end{itemize}
The detailed calculations are deferred to Section B of the Supplementary Material \textit{Proof details}.

  \begin{algorithm}[!ht]
 \caption{Graph structure learning - minimization of Equation (\ref{eq:minT})}\label{algo:minT}
 \BlankLine
   \KwIn{ $\lambda, L, \epsilon>0$.}

 {\bf Initialization:} $T_0$ the null squared $p \times p$ matrix, $k=0$ and $e=+\infty.$

  \While{$e>\epsilon$}{
  Compute $U_k=T_k - \frac{\nabla \left( \frac{1}{n} \Fro{X(I-PT_kP^T)}^2 \right)}{L}$ with Equation (\ref{eq:grad})\;
  Using Equation (\ref{eq:minTfin}), compute the current matrix $T_{k+1}=\left( (T_{k+1})_i^j \right)_{i,j}$\;
  Project $ T_{k+1}$ on $\Tri$ with Equation (\ref{eq:projT}):
  $ T_{k+1} \leftarrow \operatorname{Proj}_{\Tri}(T_{k+1})$\;
  Compute $e=\Fro{T_{k+1}-T_k}$\;
  Increase $k$: $ k\leftarrow k+1 $\;
}

  \KwOut{ $T_{k} \in \Tri$ the unique solution of (\ref{eq:minT}).}
 \end{algorithm}

\subsection{A Genetic Algorithm for a global exploration of the permutation matrices space incorporating network topologies}\label{sec:Poptim}
As the optimal $T$ can be calculated for any $P$ using Algorithm \ref{algo:minT} and with a very good approximation accuracy according to Theorem \ref{th:main}, the optimization task (\ref{eq:PT}) 
comes down to exploring the $\Perm$ space of permutation matrices in dimension $p$ and to evaluating the quality of permutation candidates $P\in \Perm$.
We first note that the number of permutation matrices is $p!$, which rules out any exact enumeration method, even for relatively small $p$. 
We propose instead to use a meta-heuristic approach, which has proven to be successful for many discrete optimization problems 
like wire-routing, transportation problems or traveling salesman problem \citep{Michalewicz96,dreo2006metaheuristics}.

Among the different meta-heuristics (Simulated annealing, Tabu search, Ant Colony,...) we focused on Genetic Algorithms (GA) because, 
despite limited convergence results \citep{Cerf98,Michalewicz96}, they were found much more efficient in problems related to ours than alternatives with more established convergence proofs (\textit{e.g.} \cite{granville1994} for simulated annealing),
while allowing the use of parallel computation. 

GAs mimic the process of natural evolution, and use a vocabulary derived from natural genetics: populations (a set of potential solutions of the optimization problem), 
individuals (a particular solution) and genes (the components of a potential solution). Each generation/iteration of the algorithm will improve the constituting elements of the population.
In short, a population made of $N$ potential solutions of the optimization problem samples the search space. 
This population is sequentially modified, with the aim of achieving a balance between exploiting the best solutions and exploring the search space, until some termination condition is met.

We use here a classical Genetic Algorithm, as described in \cite{Michalewicz96} for instance, which is based on three main operators at each iteration: selection, crossover and mutation.
The population is reduced by selection; selection shrinks the population diversity based on the individual fitness values. The crossover allows the mixing of good properties of the population to create new composite individuals. Mutations change one (or a few in more general GAs) components of the individuals to allow random space exploration. 
The complete sketch of algorithm GADAG is given in Algorithm \ref{algo:global}. A discussion on parameters to set in Algorithm \ref{algo:global} is found in Section \ref{sec-para}. The details of the different operators are given in the following.

As we show in Example \ref{ex:perm}, any $P \in \Perm$ is uniquely defined by a permutation vector of $\llbracket 1,p \rrbracket$.
Hence, we use as a the search space $\mathfrak{S}_p$ the set of permutations of $\llbracket 1,p \rrbracket$, which is a 
well-suited formulation for GAs. 

\begin{djoe}\label{ex:perm}
Consider the permutation matrix ($p=5$):
$$P=\begin{pmatrix}
0&0&0&1&0\\
0&0&0&0&1\\
0&1&0&0&0\\
0&0&1&0&0\\
1&0&0&0&0
\end{pmatrix}.$$
Then, $P$ is represented by the \begin{tikzpicture}[scale=0.2]
  \draw(0,0) rectangle (2,2);
  \draw(2,0) rectangle (4,2);
  \draw(4,0) rectangle (6,2);
  \draw(6,0) rectangle (8,2);
  \draw(8,0) rectangle (10,2);
  \draw(1,1) node{5};
  \draw(3,1) node{3};
  \draw(5,1) node{4};
  \draw(7,1) node{1};
  \draw(9,1) node{2};
  \end{tikzpicture} vector, looking at the ranks of non-null values of $P$ column by column. 
 The nodes are ranked according to their topological ordering.
\end{djoe}

Note that our problem closely resembles the classical Traveling Salesman Problem (TSP), which has been succesfully addressed by means of genetic algorithms \citep{grefenstette1985genetic,davis1991handbook}.
Identically to the TSP, we optimize over the space of permutations, which induces specific constraints for defining the crossover and mutation operators. 
However, unlike the TSP, the problem is not circular (in the TSP, the last city is connected to the first one), and the permutation here defines a hierarchy between nodes rather than a path, 
which makes the use of TSP-designed operators a potentially poor solution. As we show in the following, we carefully chose these operators in order to respect the nature of the problem at hand. In particular, we emphasize two of their desirable features: their efficiency in exploring the search space and the interpretable aspect they offer in terms of modifications on a given network or the blend of two different networks (crossover).

\paragraph{Fitness function}
Given a potential solution $p_i \in \mathfrak{S}_p$, the fitness function is defined as:
\begin{equation}\label{eq:fitness}
 J_i = J(p_i) = \frac{1}{n} \Fro{X(I - P_iT_i^*P_i^T)}^2 + \lambda \normeun{T_i^*},
 \end{equation}
with $P_i$ constructed from $p_i$ as in Example \ref{ex:perm} and $T_i^*$ the solution of Equation (\ref{eq:minT}) with $P=P_i$. 
As mentioned earlier, at each step of the proposed GA, the evaluation of the fitness function thus requires running the nested loop of our global algorithm GADAG.

\paragraph{Selection operator}
The selection operator (or survival step) consists in generating a population of N individuals from the N existing individuals by random sampling 
(with replacement, hence some individuals are duplicated and others are deleted). 
It aims at improving the average quality of the population by giving to the best potential solutions a higher probability to be copied in the intermediate population. 
We have chosen to use the classical proportional selection of \citet{Holland92}: each individuals is selected with a probability inversely proportional to its fitness value of Equation (\ref{eq:fitness}). 

\paragraph{Crossover operator}
A crossover operator generates a new set of potential solutions (children) from existing solutions (parents). Crossover aims at achieving at the same time (i) a good exploration 
of the search space by mixing the characteristics of the parents to create potentially new ones while (ii) preserving some of the desirable characteristics of the parents. By desirable features, we mean features of the network which lead to good fitness values, and which in turn are favored by selection over the generations.
The crossover population (set of parents) is obtained by selecting each individual of the population with a probability $p_{xo}$; the parents are then paired randomly.  

We have chosen the \textit{order-based} crossover, originally proposed for the TSP \citep[Chapter 10]{Michalewicz96}, which is defined as follows. 
Given two parents $p_1$ and $p_2$, a random set of crossover points are selected, which we denote $\Omega$. It consists in a permutation of $k$ elements taken from $\llbracket 1,p \rrbracket$, with $k$ uniformly drawn between 0 and $p$. A first child $C_1$ between $p_1$ and $p_2$ is then generated by:
\begin{enumerate}
\item fixing the crossover points of $p_1$,
\item completing $C_1$ with the missing numbers in the order they appear in $p_2$.
\end{enumerate}

\begin{djoe}\label{ex:cross}
Consider the two following parents:

\begin{center}
 \begin{tikzpicture}[scale=0.3]
 \draw(-2,1) node{$p_1$};
  \draw(0,0) rectangle (2,2);
  \draw(2,0) rectangle (4,2);
  \draw(4,0) rectangle (6,2);
  \draw(6,0) rectangle (8,2);
  \draw(8,0) rectangle (10,2);
  \draw(10,0) rectangle (12,2);
  \draw(10,0) rectangle (14,2);
  \draw(10,0) rectangle (16,2);
  \draw(10,0) rectangle (18,2);
  \draw(10,0) rectangle (20,2);
  \draw(1,1) node{\textbf{\textcolor{red}{4}}};
  \draw(3,1) node{3};
  \draw(5,1) node{10};
  \draw(7,1) node{7};
  \draw(9,1) node{5};
  \draw(11,1) node{\textbf{\textcolor{red}{9}}};
    \draw(13,1) node{1};
  \draw(15,1) node{\textbf{\textcolor{red}{2}}};
  \draw(17,1) node{6};
 \draw(19,1) node{\textbf{\textcolor{red}{8}}};

  \draw(-2,-3) node{$p_2$};
  \draw(0,-4) rectangle (2,-2);
  \draw(2,-4) rectangle (4,-2);
  \draw(4,-4) rectangle (6,-2);
  \draw(6,-4) rectangle (8,-2);
  \draw(8,-4) rectangle (10,-2);
  \draw(10,-4) rectangle (12,-2);
  \draw(10,-4) rectangle (14,-2);
  \draw(10,-4) rectangle (16,-2);
  \draw(10,-4) rectangle (18,-2);
  \draw(10,-4) rectangle (20,-2);
  \draw(1,-3) node{6};
  \draw(3,-3) node{1};
  \draw(5,-3) node{\textbf{\textcolor{red}{9}}};
  \draw(7,-3) node{\textbf{\textcolor{red}{4}}};
  \draw(9,-3) node{10};
  \draw(11,-3) node{\textbf{\textcolor{red}{2}}};
    \draw(13,-3) node{\textbf{\textcolor{red}{8}}};
  \draw(15,-3) node{3};
  \draw(17,-3) node{7};
 \draw(19,-3) node{5};
  \end{tikzpicture}
  \end{center}

Assume that the crossover points randomly chosen are $4$, $9$, $2$ and $8$ (in bold red above). Then, the child $C_1$ is defined by 
inheriting those points from $p_1$ and filling the other points in the order they appear in $p_2$:
   \begin{center}
 \begin{tikzpicture}[scale=0.3]
 \draw(-2,-3) node{$C_1$};
  \draw(0,-4) rectangle (2,-2);
  \draw(2,-4) rectangle (4,-2);
  \draw(4,-4) rectangle (6,-2);
  \draw(6,-4) rectangle (8,-2);
  \draw(8,-4) rectangle (10,-2);
  \draw(10,-4) rectangle (12,-2);
  \draw(10,-4) rectangle (14,-2);
  \draw(10,-4) rectangle (16,-2);
  \draw(10,-4) rectangle (18,-2);
  \draw(10,-4) rectangle (20,-2);
  \draw(1,-3) node{\textbf{\textcolor{red}{4}}};
  \draw(3,-3.3) node{*};
  \draw(5,-3.3) node{*};
  \draw(7,-3.3) node{*};
  \draw(9,-3.3) node{*};
  \draw(11,-3) node{\textbf{\textcolor{red}{9}}};
  \draw(13,-3.3) node{*};
  \draw(15,-3) node{\textbf{\textcolor{red}{2}}};
  \draw(17,-3.3) node{*};
 \draw(19,-3) node{\textbf{\textcolor{red}{8}}};

   \draw(-2,1) node{$p_2$};
  \draw(0,0) rectangle (2,2);
  \draw(2,0) rectangle (4,2);
  \draw(4,0) rectangle (6,2);
  \draw(6,0) rectangle (8,2);
  \draw(8,0) rectangle (10,2);
  \draw(10,0) rectangle (12,2);
  \draw(10,0) rectangle (14,2);
  \draw(10,0) rectangle (16,2);
  \draw(10,0) rectangle (18,2);
  \draw(10,0) rectangle (20,2);
  \draw(1,1) node{6};
  \draw(3,1) node{1};
  \draw(5,1) node{9};
  \draw(7,1) node{4};
  \draw(9,1) node{10};
  \draw(11,1) node{2};
    \draw(13,1) node{8};
  \draw(15,1) node{3};
  \draw(17,1) node{7};
 \draw(19,1) node{5};
 \draw(5,1) node{/};
 \draw(7,1) node{/};
 \draw(11,1) node{/};
 \draw(13,1) node{/};
 \draw[->,thick] (1,0)--(3,-2);
 \draw[->,thick] (9,0)--(7,-2);
 \draw[->,thick] (15,0)--(9,-2);
 \draw[->,thick] (17,0)--(13,-2);
 \draw[->,thick] (19,0)--(17,-2);
 \draw[->,thick] (3,0)--(5,-2);

 \draw(9,-6) node{$\Downarrow$};
  \draw(0,-8) rectangle (2,-10);
  \draw(2,-8) rectangle (4,-10);
  \draw(4,-8) rectangle (6,-10);
  \draw(6,-8) rectangle (8,-10);
  \draw(8,-8) rectangle (10,-10);
  \draw(10,-8) rectangle (12,-10);
  \draw(10,-8) rectangle (14,-10);
  \draw(10,-8) rectangle (16,-10);
  \draw(10,-8) rectangle (18,-10);
  \draw(10,-8) rectangle (20,-10);
  \draw(1,-9) node{4};
  \draw(3,-9) node{6};
  \draw(5,-9) node{1};
  \draw(7,-9) node{10};
  \draw(9,-9) node{3};
  \draw(11,-9) node{9};
      \draw(13,-9) node{7};
  \draw(15,-9) node{2};
  \draw(17,-9) node{5};
 \draw(19,-9) node{8};

   \end{tikzpicture}
  \end{center}

\end{djoe}

\begin{figure*}[!ht]
\begin{center}
\includegraphics[trim= 7.5cm 2.5cm 6.5cm 2cm, clip, width=5cm]{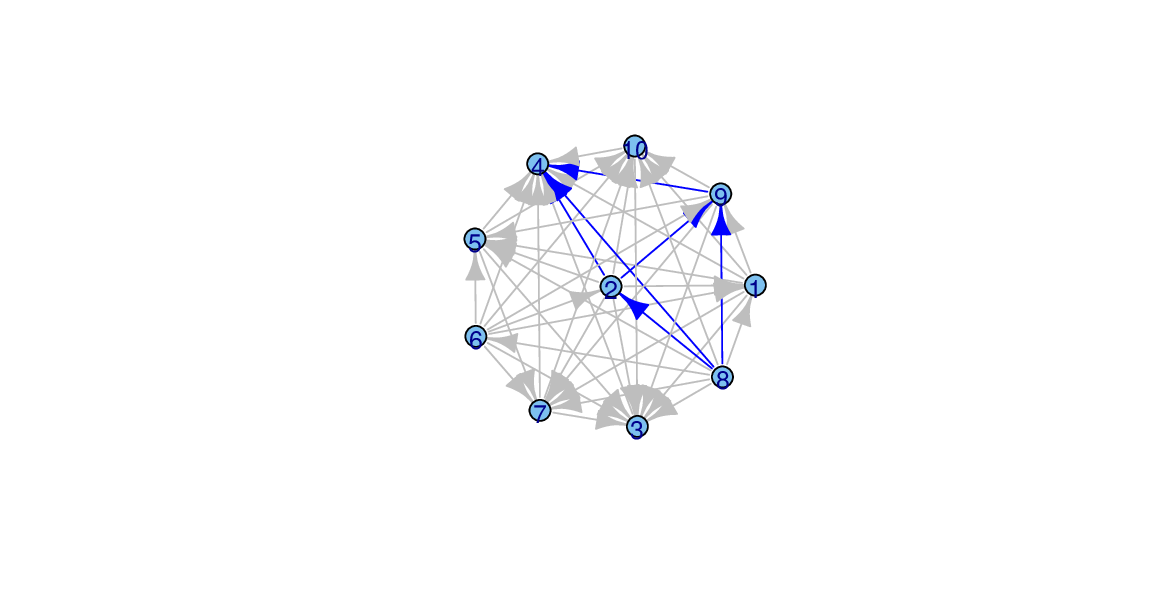}
\includegraphics[trim= 7.5cm 2.5cm 6.5cm 2cm, clip, width=5cm]{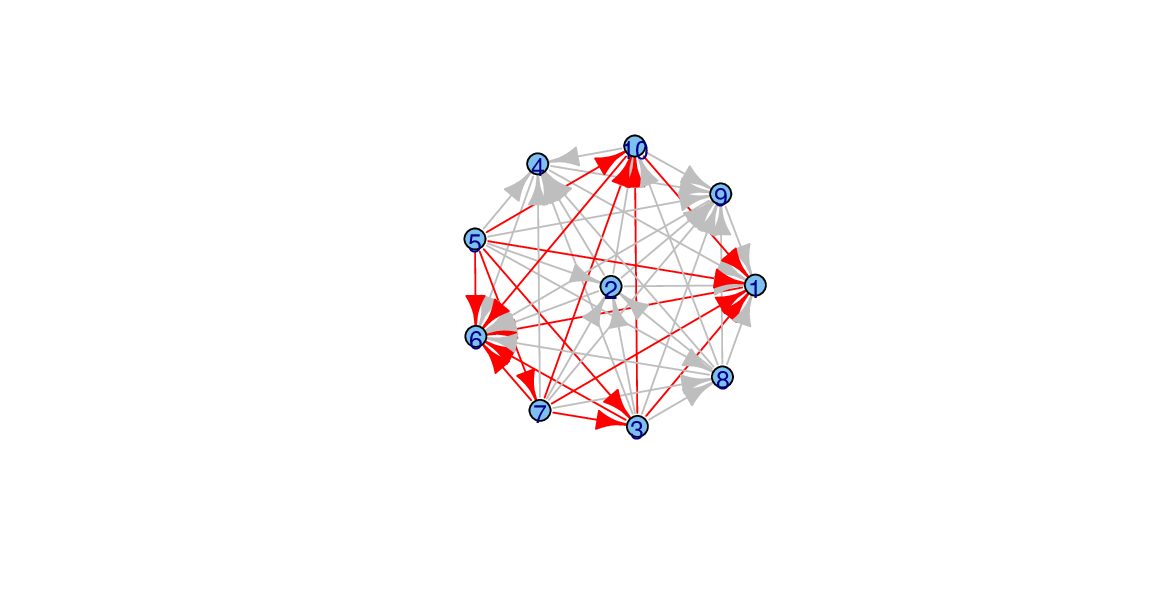}
\includegraphics[trim= 7.5cm 2.5cm 6.5cm 2cm, clip, width=5cm]{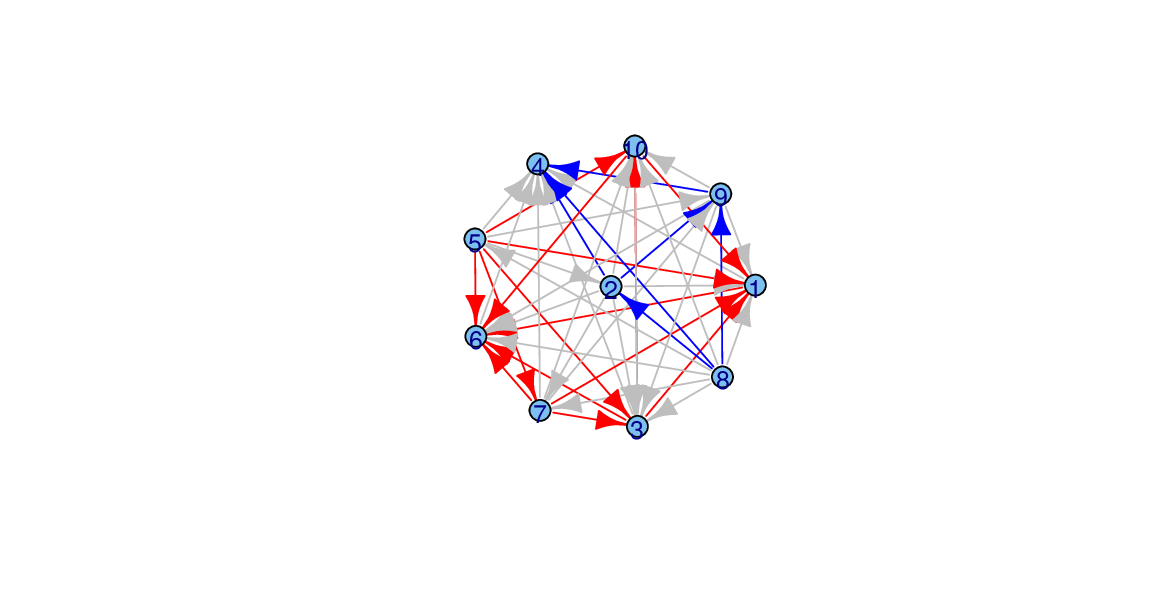}
\end{center}
\caption{Graphical representation of crossover between two 10-node graphs. The two parental graphs are represented on the left. The third graph, on the right, is obtained by combining the blue and red part of its parents using the crossover operator.
}\label{fig:croisement10}
\end{figure*}

From a graphical point of view, a crossover between $p_1$ and $p_2$, which encode two complete graphs $\mathcal{G}_{P_{1}} $ and $\mathcal{G}_{P_{2}} $, 
constructs two new graphs. One of them, $\mathcal{G}_{C_{1}} $ is composed of the sub-graph of $\mathcal{G}_{P_{1}} $ induced by the  set of crossover points $\Omega$ and the sub-graph of $\mathcal{G}_{P_{1}} $ induced by 
the complementary set $\Omega^C$ of $\Omega$ in $\llbracket 1,p \rrbracket$ (see Figure \ref{fig:croisement10}). The second child graph $\mathcal{G}_{C_{2}} $ is obtained in an identical manner by reversing the roles played by the two parental graphs.

\paragraph{Mutation operator}
Mutation operators usually correspond to the smallest possible change in an individual (unary operator). We thus define it as an alteration of
two neighbouring genes (see Example \ref{ex:mutation}). Graphically, a mutation consists in switching the arrowhead of an edge between two nodes.
Mutation is applied to each child with probability $p_m$.

\begin{djoe}\label{ex:mutation}
A possible mutation for the first child of Example \ref{ex:cross} is to swap the genes ``1'' and ``10'' (in bold red below):
\begin{center}
 \begin{tikzpicture}[scale=0.3]
 \draw(-2,1) node{$M_1$};
  \draw(0,0) rectangle (2,2);
  \draw(2,0) rectangle (4,2);
  \draw(4,0) rectangle (6,2);
  \draw(6,0) rectangle (8,2);
  \draw(8,0) rectangle (10,2);
  \draw(10,0) rectangle (12,2);
  \draw(10,0) rectangle (14,2);
  \draw(10,0) rectangle (16,2);
  \draw(10,0) rectangle (18,2);
  \draw(10,0) rectangle (20,2);
  \draw(1,1) node{4};
  \draw(3,1) node{6};
  \draw(5,1) node{\textbf{\textcolor{red}{1}}};
  \draw(7,1) node{\textbf{\textcolor{red}{10}}};
  \draw(9,1) node{3};
  \draw(11,1) node{9};
    \draw(13,1) node{7};
  \draw(15,1) node{2};
  \draw(17,1) node{5};
 \draw(19,1) node{8};

  \draw[<-,thick] (5,0)--(5,-1); 
 \draw (5,-1) --(7,-1);
 \draw[->,thick] (7,-1) -- (7,0);
 \end{tikzpicture}
 \end{center}
\end{djoe}

\paragraph{Stopping criterion}
Two quantities are monitored along the iterations: the heterogeneity of the population and the value of the objective function. 

For the first indicator, we use the Shannon entropy, defined for each rank position $j\in \llbracket 1,p \rrbracket$ as:
$$H_j = -\sum_{i=1}^p \frac{N_{i,j}}{N} \log \left( \frac{N_{i,j}}{N} \right),$$
where $N_{i,j}$ is the number of times when $i$ appears in position $j$. 
$H_j=0$ if all the individuals ``agree'' on the position of a node and the population is perfectly homogeneous at this node. On the contrary, it is maximum when we observe a uniform distribution of the different nodes at a given position and the population is in this case at a maximum of heterogeneity or disorder for this position.
The algorithm stops if the population entropy value $ H = \sum_{j=1}^N H_j $ drops below a threshold since $ H=0 $ if all the individuals are identical. 
A second criterion can terminate GADAG if difference in the average fitness (denoted $ \bar{J} $ thereafter) of the population between a given number of consecutive iterations, does not change by more than a predefined threshold.

  \begin{algorithm}[!ht]
 \caption{GADAG overview}\label{algo:global}
 \BlankLine
   \KwIn{ $p_{xo}, p_m, \epsilon_H>0, \epsilon_J>0$, $k_{max}>0$, $i_{max}>0$, $\lambda$, $L$.}

 {\bf Initialization:} Generate the initial population $\mathcal{P}_0$ with $N$ permutations of $\llbracket 1,p \rrbracket$, $k=0$ and $e_J=+\infty$.

  \While{$H>\epsilon_H \quad \& \quad e_J>\epsilon_J \& \quad k\leq k_{max}$}{
  Generate $\mathcal{P}_{k+1}$ as a random \textbf{selection} of $N$ individuals from $\mathcal{P}_k$\;
  Pick an even subset $\mathcal{P}_{xo}$ of $\mathcal{P}_{k+1}$ (each individual of $\mathcal{P}_{k+1}$ selected with probability $p_{xo}$)\;
  Perform \textbf{crossover} on $\mathcal{P}_{xo}$ by randomly pairing the individuals\;
  \textbf{Mutate} each obtained individual with probability $p_m$ \;
  \textbf{Evaluate} the new individuals $\mathcal{P}_m$ by running Algorithm \ref{algo:minT}\;
  Replace $\mathcal{P}_{xo}$ by $\mathcal{P}_m$ in $\mathcal{P}_{k+1}$\;
  Compute the Shannon entropy $H$ and the difference in the average fitness $e_J=\max_{0\leq i\leq i_{max} } \left( \bar{J}(\mathcal{P}_{k+1}) - \bar{J}(\mathcal{P}_{k-i})\right)$\;
  Increase k: $k\leftarrow k+1$\;
}
 \end{algorithm}

\section{Numerical experiments}\label{sec:experiments}
This section is dedicated to experimental studies to assess practical performances of our method through two kinds of datasets. 
In a first phase, the aim of these applications is to show that GADAG has a sound behavior on simulated toy data with a variety of different settings. 
In a second phase, we demonstrate the ability of our algorithm to analyse datasets that
mimic the activity of a complex biological system,
and we compare it to other state-of-the art methods. The competing methods are presented in Section \ref{sec-comp}.
In Section \ref{sec-para}, we present the calibration of the Genetic Algorithm parameters. Section \ref{sec-perf} introduces the measures we used to assess the merits of the methods. 
 Experimental results are then detailed in  
Section \ref{sec-simus1} for the simulated high-dimensional toy datasets and in 
Section \ref{sec-simus2} for the dataset with features encountered in real situations.

All experiments have been performed on \texttt{R} \citep{Rproject} using the package \texttt{GADAG} \citep{champion2017}.
The computational times reported in Section \ref{sec-simus1} correspond to a Windows 7 laptop computer with 
8 threads on a 4-core hyperthreaded 2.50GHz processor, with 4GB of RAM.

 \subsection{Algorithm parameters}\label{sec-para}
Running the procedure of Algorithm \ref{algo:global} requires to define parameters of the outer loop, which generates our population of $P$'s, and of the nested loop to find the optimal $T^*$. 
The evaluation of the Lipschitz gradient constant $L$, used to find the optimal graph structure $T^*$, is known as a hard established problem in optimization. 
Some authors propose to choose an estimate of $L$ from a set of possible values \citep{Jones93,L06}, to estimate local Lipschitz constants \citep{L95}, or to set it a priori to a fixed value \citep{L09,Horst95}.
Here, observing Equation (\ref{eq:grad}), a major bound for $L$ is given by:
$$L\leq \frac{2}{n}\Fro{X^TX}.$$
We found that setting $L$ to this bound worked well in practice in all our scenarios.

Five parameters need to be tuned to run the Genetic Algorithm:
the crossover rate $ p_{xo} $, the mutation rate $ p_m $, the constant of the stopping criteria $ \epsilon_{H}$ and $ \epsilon_{J} $ and the size of the population $ N$. For the first four parameters, we observed that their value had a limited effect on the efficiency, hence we chose commonly used values in the literature (see Table~\ref{tab:param}). 
The size of the population has a more complex effect and has been investigated in several prospective papers (\textit{e.g.} \citealt{Schaffer89, alander1992, piszcz2006, ridge2007}) but without providing a definitive answer to the problem. 
In our simulation study, we chose as a rule-of thumb $N=5p$, which was found as a good compromise between computational cost and space exploration on several experiments.

The complete parameter settings used in our experiments are reported in Table~\ref{tab:param}.

\begin{table}[!ht]
\begin{center}
 \caption{Algorithm parameter settings}\label{tab:param}
 \begin{tabular}{l l}
 \hline\noalign{\smallskip}
Parameter & Value \\
\noalign{\smallskip}\hline\noalign{\smallskip}
$p_{xo}$ & 0.25\\
 $p_{m}$ & 0.5\\
 $N$ & $5 \times p$\\
 $L$ & $\frac{2}{n}\Fro{X^TX}$\\
 max. nb. of eval. &$10^4$ \\
 $\epsilon_{H}$ & $10^{-6}$\\
  $ \epsilon_{J} $&$10^{-4}$ \\
\noalign{\smallskip}\hline
 \end{tabular}
 \end{center}
\end{table}

\subsection{Performance metrics}\label{sec-perf}
A classical performance measure for graph inference methods consists in comparing predicted interactions with the known edges in the true graph $\mathcal{G}_0$ using precision versus recall (P/R) curves.
We denote by TP, FP, FN and TN, the true positive (correctly predicted) edges, the false positive (inferred by mistake) edges, the false negative (missed) edges, and the true negative (correctly non-predicted)
edges. 
The recall, defined as $\frac{TP}{TP+FN}$, measures the
power (or sensitivity) of reconstruction of non-zero elements of the true matrix $G$ (or equivalently of the true network) for one method, whereas the precision, equal to $\frac{TP}{TP+FP}$, measures the
accuracy of the reconstruction. The closer to one the precision and the recall the better. 

P/R curves represent the evolution of those quantities when varying the sparsity of the methods.
GADAG
is based on penalized optimization: 
it seeks linear dependencies between the variables with a controlled level of parsimony ($\lambda$ in Equation (\ref{eq:lv2})).
For $\lambda$ varying from $0$ (complete graph) to $+\infty$ (empty graph), it thus produces a list of edges successively introduced in the model.
This list of edges defines the precision versus recall curve.
As a summary performance measurement, we also computed the classical area under the P/R curve (AUPR).

\subsection{Exploratory analysis on toy datasets} \label{sec-simus1}
We first considered simulated data from randomly generated networks with different characteristics
in order to assess the capabilities and limits of our algorithm. Given a number of nodes $p$,
a random set of $s$ edges were generated, and the non-zero parameters of the matrix $G_0$ 
associated to the corresponding DAG were uniformly sampled between 0 and 1. 
Using this graph, we generated $N$ observations following the hypotheses of 
Gaussian, homoscedastic and centred error. We then ran GADAG on this dataset to recover the graph. Note that other assumptions presented in Section \ref{sec-hyp} may not be fulfilled here, but we aimed at evaluating the robustness of GADAG for recovering DAGs in such scenario.

In our experiments, we varied the number of nodes $p$, of edges $s$ and of available observations $n$.
We chose four different settings $p=50, 100, 500$ and $1,000$ with $n/p$ varying from $100\%$ to $10\%$ and $s/p$ from $100\%$ to $400\%$.
Unless otherwise stated, all experiments were replicated $50$ times each 
and results were averaged over these replicates. Averaged computational times correspond to one iteration of GADAG, for a fixed parameter of penalization $\lambda$.

\begin{table*}[!ht]
\caption{Performance results of GADAG on a toy dataset with different characteristics (number of nodes $p$, number of edges $s$ and sample size $n$) in terms of area under the Precision vs. Recall curve (a) and computational time, in seconds (b).
All results are averaged over 50 replicates ($^*$ 5 replicates only as the running time was 1 day per network).
}\label{tab:diff}

(a)

\begin{center}
\begin{small}
 \begin{tabularx}{0.95\textwidth}{c c c c c c c c c c c c c c}
\toprule
& & \multicolumn{3}{c}{p=50} & \multicolumn{3}{c}{p=100} & \multicolumn{3}{c}{p=500} &\multicolumn{3}{c}{p=1,000$^*$}  
\\
 \cmidrule(lr){3-5} \cmidrule(lr){6-8} \cmidrule(lr){9-11}\cmidrule(lr){12-14}
& $s/p$ & $100\%$ & $200\%$ & $400\%$ &$100\%$ & $200\%$ & $400\%$& $100\%$ & $200\%$ & $400\%$ &  $100\%$ & $200\%$ & $400\%$\\
\midrule
& $100\%$& 0.72 & 0.65 & 0.65 & 0.69 & 0.63& 0.64 & 0.79&0.79 & 0.70& 0.86& 0.80 & 0.65\\
$n/p$ & $50\%$ & 0.53& 0.46 & 0.49&0.56 & 0.53 & 0.53&0.75&0.76&0.66&0.82 &0.75&0.63\\
& $10\%$ & 0.02 & 0.03 & 0.05 & 0.04 & 0.06 & 0.09 & 0.51&0.53&0.41&0.65 & 0.60&0.49\\
\bottomrule
 \end{tabularx}
 \end{small}
\end{center}

(b)

\begin{center}
\begin{small}
 \begin{tabularx}{0.95\textwidth}{c c c c c c c c c c c c c c}
\toprule
& & \multicolumn{3}{c}{p=50} & \multicolumn{3}{c}{p=100} & \multicolumn{3}{c}{p=500} &\multicolumn{3}{c}{p=1,000$^*$}  
\\
 \cmidrule(lr){3-5} \cmidrule(lr){6-8} \cmidrule(lr){9-11}\cmidrule(lr){12-14}
& $s/p$ & $100\%$ & $200\%$ & $400\%$ &$100\%$ & $200\%$ & $400\%$& $100\%$ & $200\%$ & $400\%$ &  $100\%$ & $200\%$ & $400\%$\\
\midrule
& $100\%$& 0.55 & 0.43 & 0.41& 8.69 & 9.87&7.51& 253&258&214&1,640&1,500&1,230\\
$n/p$ & $50\%$ & 0.48 & 0.48 & 0.41 & 8.63& 8.90& 6.62&250&258&172&1,520 & 1,590 &1,270 \\
& $10\%$ &0.31 & 0.30 & 0.36 & 6.52 & 5.81 & 5.41&183&183&165&1,590 & 1,550 &1,290\\
\bottomrule
 \end{tabularx}
 \end{small}
\end{center}

\end{table*}

Results, in terms of area under the P/R curves and computational time are sumarized in Table \ref{tab:diff}.
We can first remark a crude decrease of performance results when the number of samples is very small ($p=50$ and $100$, $n/p=10\%$, so respectively $5$ and $10$ samples).
In that case, GADAG is incapable of recovering any signal (AUPR $<10\%$). 
When the sample size is of the order of $p$ ($n/p=100\%$, first row of Table \ref{tab:diff} a), GADAG works well, although it is clearly a favorable case, far from the high-dimensional one. 
With half of the samples ($n/p=50\%$), performance remains satisfactory (AUPR around $50\%$, or more), which is critical since 
this situation corresponds to realistic biological studies, where subsets of genes (\textit{i.e.} nodes) are preselected beforehand.
Interestingly, $p=500$ and $1,000$ work better than smaller values of $p$ since the number of samples is larger to estimate the graph that generated the data. Indeed, for a given number of samples, \textit{e.g.}
$n=50$, performance results slightly decrease from $65\%$ ($p=50$) and $55\%$ ($p=100$) to $50\%$ ($p=500$). GADAG is thus not considerably affected by a relative increase of dimensionality. For large graphs ($p=500$ and $1,000$), even if $n/p\leq 10\%$, it suceeds in recovering them, which makes it a competitive algorithm with regards to other state-of-the-art approaches.


An interesting remark is that the number of edges $s$ does not significantly  change numerical results (see each row of Table \ref{tab:diff} a),
although GADAG succeeds slightly better in estimating sparser graphs. This may be due to the particular structure of our algorithm, 
which looks for topological ordering between nodes (genetic algorithm) and then makes the inferred graph sparse.

Concerning computational time, we can finally note that growing the dimension $p$ clearly makes the problem harder to solve : each call to GADAG requires more than $300s$ for hundred of nodes and $1,500s$ for thousand of nodes.

\subsection{DREAM data analysis}
The second type of datasets we used mimic activations and regulations that occur in gene regulatory networks. It is provided by the DREAM4 challenge on ``In Silico Network Challenge''. 
Note that although plausibly simulated, DREAM4 data sets are not real biological data sets. However, the used network structures (five in total) were extracted from {\it E. coli} and {\it S. cerevisae} 
-two biological model organisms- trancriptional networks. These networks contain cycles, but self-loops were discarded. The gene expression observations were not simulated by an equal noise Gaussian multivariate model, stochastic differential equations were used to mimic the kinetic laws of intricate and intertwined gene regulations. 
In addition to the biological noise simulated from the stochastic differential equations, technical noises were added to reproduce actual gene measurement noise. 
All data sets were generated by the {\tt GNW} software \citep{marbach2009}.

Working with simulated networks, we are able to quantitatively and objectively assess the merit of competing methods in terms of true predictions (true positive TP and true negative TN) vs. incorrect predictions (false positive FP and false negative FN)
edges. 
While the analysis of a real data set is certainly the final goal of a methodology motivated by a real problem like ours, there are only imprecise ways of validating a method when analysing a real data set. 
Well known systems are often small and even if knowledge has accumulated on them, these can be noisy and difficult to gather to obtain a fair picture of what can adequately be considered as sets of true positive and true negative sets of edges. 
Even if the data generation process of the DREAM4 In Silico Network Challenge is completely understood, no existing method is able to predict all regulatory relationships, but at the price of including many false positive predictions. 
The DREAM4 datasets we considered have $ p = 100 $ nodes and only $ n = 100 $ observations 
making it a a very challenging task. 

\subsubsection{Comparison to state-of-the art}\label{sec-comp}
We compare GADAG to five state-of-the-art inference methods. 
Among them, the Genie3 method \citep{Huynh10}, based on random forests, was the best performer of one of the DREAM4 sub-challenges, while the BootLasso \citep{Allouche13} was one of the key components of the best performing approach of one of the DREAM5 sub-challenges \citep{Allouche13}. The two methods decompose the prediction of the network into $ p $ feature selection sub-problems. In each of the $ p $ sub-problems, one of the node is predicted from the other ones using random forests \citep{Breiman01} for Genie3 or a bootstrapped version of the Lasso \citep{Bach08} for BootLasso.
For the random forest approach, 
parents of each node were detected as most significant explanatory variables according to a variance reduction criterion in a regression tree framework. The process was repeated on a randomized set of trees, which made up the so-called random forest. 
This method allowed us to derive a ranking of the importance of all variables for the target by averaging the scores over all the trees.
We used the R package {\tt randomforest} \citep{liaw2002} for our results.
The Lasso is a $\ell_1$-norm penalization technique for solving linear regression. Following the works of \citet{Bach08}, BootLasso uses bootstrapped estimates of the active regression set based on a Lasso penalty: only those variables that are selected in every bootstrap are kept in the model. In both cases, actual coefficient values are estimated from a straightforward least square procedure. Note that we slightly relax the condition for a variable to be included in the model, a variable
was selected at a given penalty level if more than $ 80\% $ of bootstrapped samples led to selecting it in the model \citep{Allouche13}. 

We also compare our algorithm to three classical methods for Bayesian Networks (BNs) modelling. BNs are graphical models \citep{Pearl09} defined by a DAG and parameters that set quantitative relationships between variables.
Algorithms devoted to structure and parameter learning in BNs either aim at maximising a score that reflects the fit of the data to the learnt structure, or test for independencies between variables. They are often used as references in a gene regulatory network inference context \citep{Tsamardinos06}, athough mainly for moderate size networks. The first compared algorithm we used is the PC-algorithm \citep{Spirtes00}, a popular constraint-based method that drastically reduces the number of conditional independence tests. It first builds the skeleton of the graph by removing edges from a complete undirected graph before determining the orientation of the edges, when possible. A large number of implementations of the PC-algorithm exists. The numerical results presented here were obtained using the pcAlg function of the R-package {\tt pcalg}, based on standard correlation estimates for conditional indepence testing.
We also ran ARACNE \citep{Margolin06}, an improved version of minimum-weight spanning tree that uses the information inequality to eliminate the majority of indirect relationships between variables. We used the ARACNE function of
the R-package {\tt bnlearn}.
We finally compare GADAG to the Greedy Equivalence Search (GSE) algorithm  \citep{Chickering02}, implemented in the R-package {\tt pacalg}, which heuristically searches in the space of equivalent classes the model with the highest Bayesian score.

To compare our algorithm with these competing methods, we used the P/R curves presented in Section \ref{sec-perf}. As GADAG, BootLasso leads to a sparse inferred graph while controlling the level of parsimony, which builds the P/R curve.
Genie3 produces as an output a ranked list of regulatory interactions, which corresponds to the edges of the inferred graph. Edges are the successively introduced with decreasing confidence scores to produce the random forest P/R curve. For the PC and the GSE algorithms, inherent parameters regulating the sparsity of the produced graphs helped us to define such curves. 
Note that the implementation we used for running ARACNE was only able to produce a final network prediction (interaction ranking is not available).

\subsubsection{Numerical results} \label{sec-simus2}
\begin{figure*}[!ht]
 \centering
 \begin{tabular}{cc}
 A & B \\
\includegraphics[trim= 1mm 5mm 10mm 20mm, clip, width=.45\textwidth]{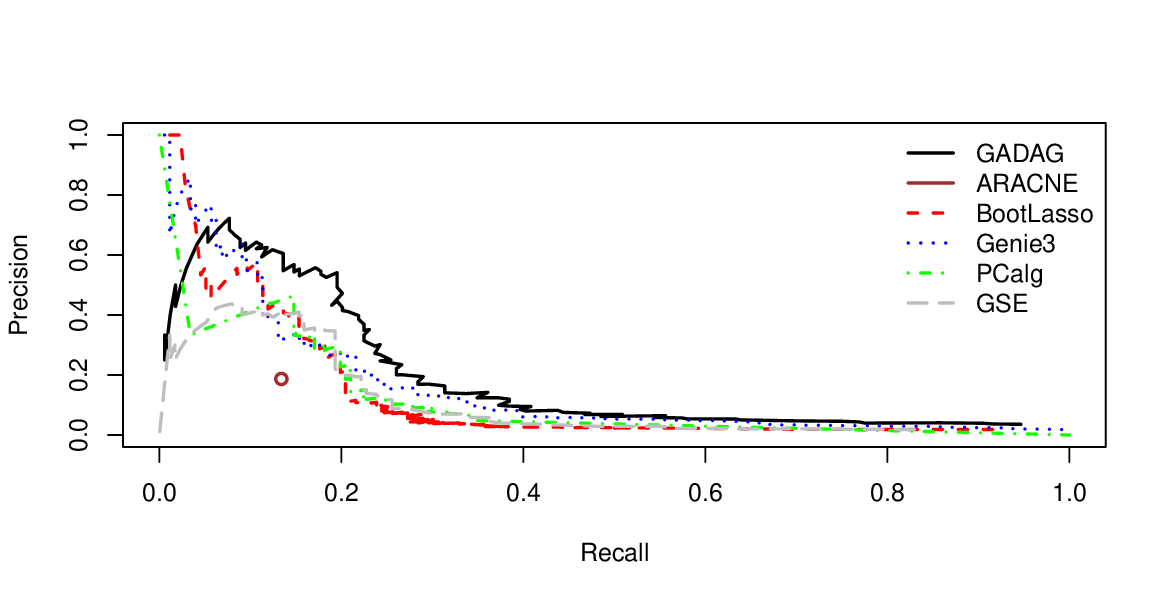} &
\includegraphics[trim= 1mm 5mm 10mm 20mm, clip, width=.45\textwidth]{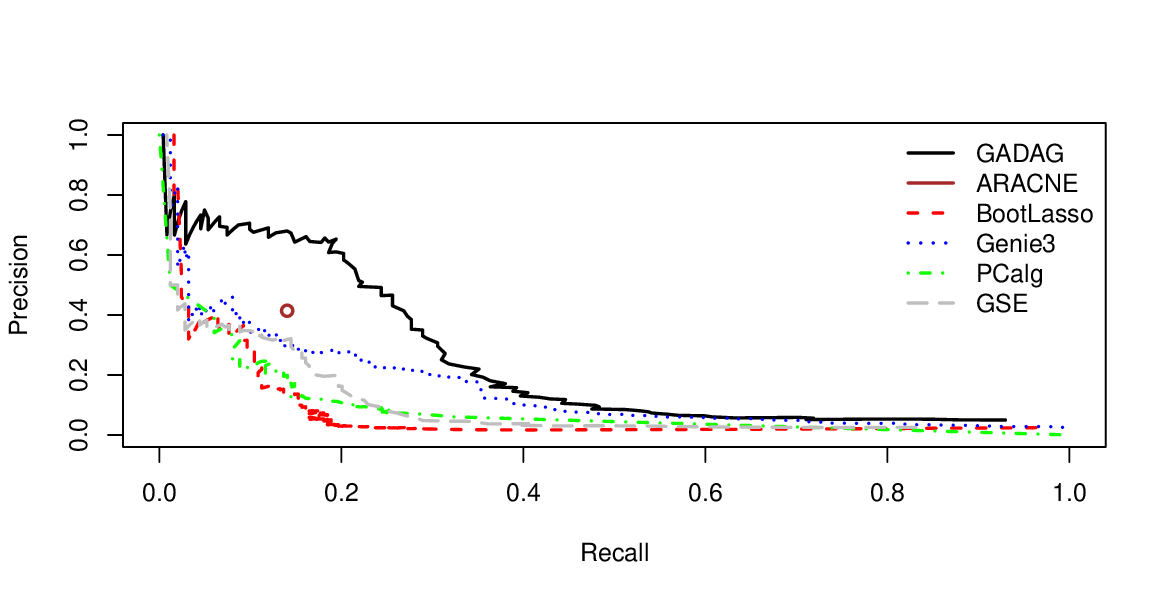}  \\
C & D \\
\includegraphics[trim= 1mm 5mm 10mm 20mm, clip, width=.45\textwidth]{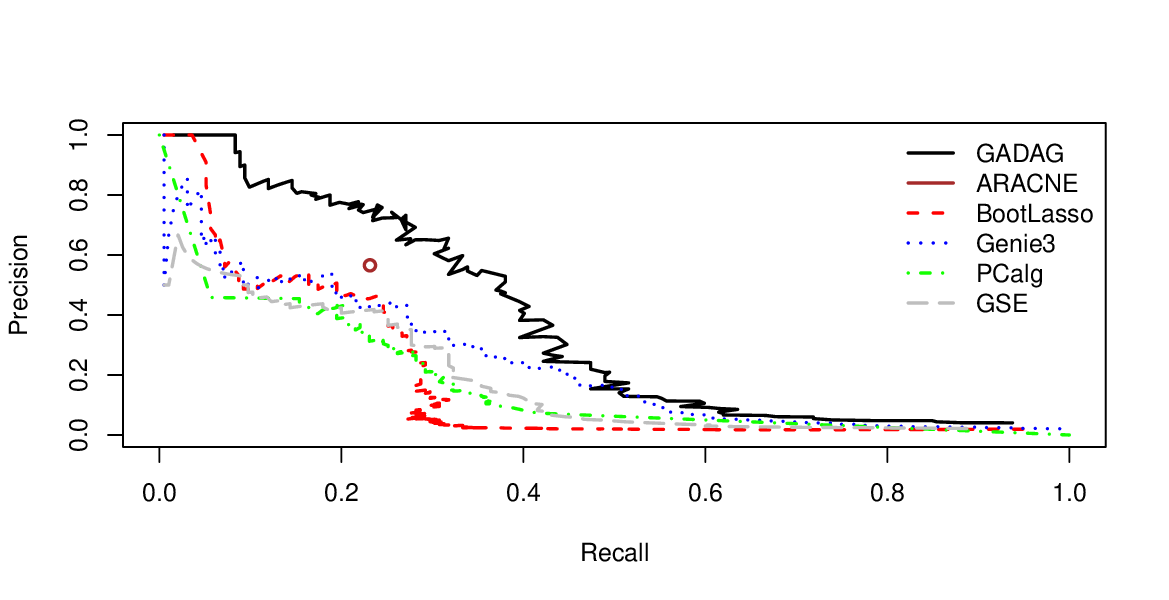} &
\includegraphics[trim= 1mm 5mm 10mm 20mm, clip, width=.45\textwidth]{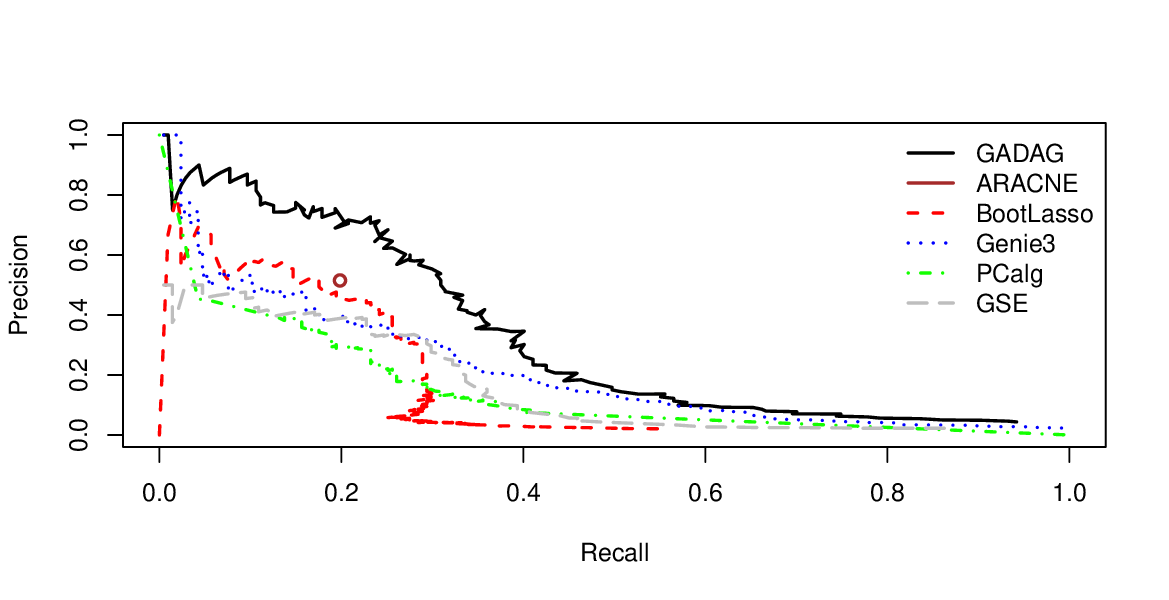} \\
E & \\
\includegraphics[trim= 1mm 5mm 10mm 20mm, clip, width=.45\textwidth]{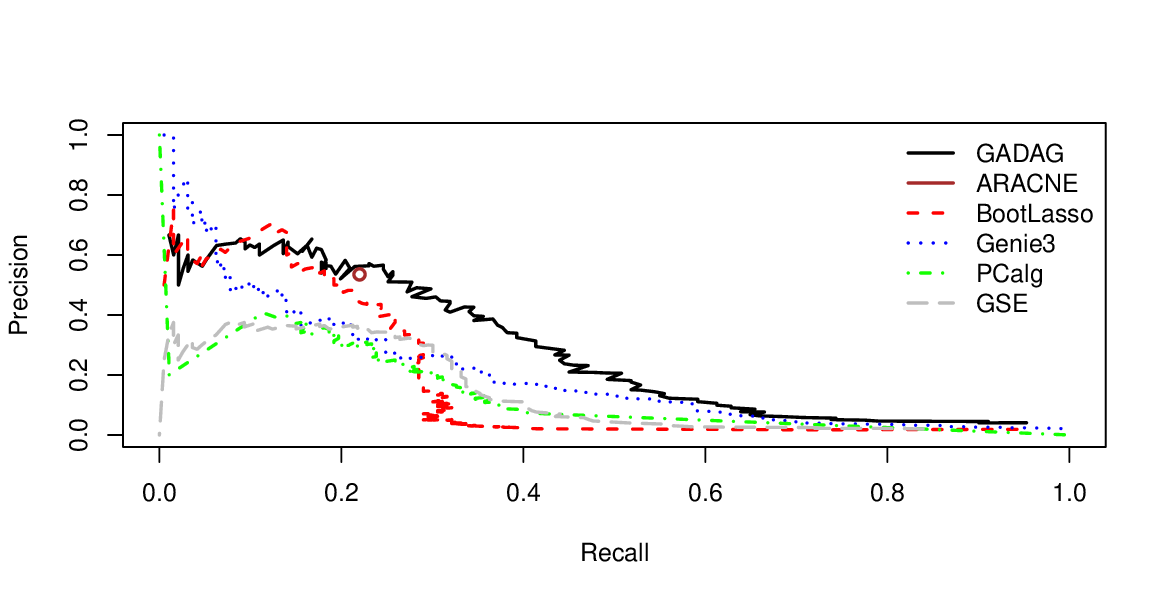}
 \end{tabular}
\caption{P/R curves for the five Dream networks and the five compared methods.
}
\label{fig:res1}
\end{figure*}

The P/R curves for the five DREAM problems are shown in Figure \ref{fig:res1}. Each curve corresponds to one of the five networks used in the challenge. In general, for all the problems the five methods are able to achieve a precision equal to one (that is, to include only true edges), but these correspond to overly sparse graphs (very small recall). Conversely, a recall equal to $ 1 $ can only be reached by adding a large number of FP edges, whatever the method we consider, even if some fail earlier than others. The main differences between the methods appear on the leftmost part of the P/R curves, especially those of Figure \ref{fig:res1} B, C and D: while the precision of BootLasso, Genie3, PCalg and GSE drops rapidly with a slow increases in recall above $ 20 \% $ recall, it remains higher for GADAG. Hence, its first predicted edges are at least as accurate than those of the four other methods and it produces a larger set of reliable edges.
For graphs of lesser sparsity, none of the five methods is really able to identify clearly reliable edges. Large number of FP edges are produced to achieve a recall higher than $ 60 \% $.





For Networks 1 and to a lesser extent 5 (Figure \ref{fig:res1} A and E), GADAG recovers with more difficulty the first true edges than other methods, with a high level of FP edges at the beginning of the curve (low precision and low recall).
However, as soon as the recall exceeds the 10\%, \textit{resp.} 15\%, for graph A, \textit{resp.} for graph E, GADAG performance is again superior to that other methods.



Table \ref{table:AUPR} gives the areas under the P/R curves for all methods and networks. For this indicator, GA significantly outperforms  the state-of-the-art methods for all networks. 

%

\begin{table}[!ht]
 \caption{Area under the Precision vs. Recall curve for all networks and methods (except ARACNE).}
 \label{table:AUPR}
\begin{center}
 \begin{tabular}{l l l l l l}
 \hline \noalign{\smallskip}
Method & Net 1 & Net 2 & Net 3 & Net 4 & Net 5 \\
\noalign{\smallskip}\hline\noalign{\smallskip}
  GADAG& \textbf{0.182}	& \textbf{0.236} & \textbf{0.348} & \textbf{0.317} & \textbf{0.267}\\
  Genie3 & 0.154 & 0.155 & 0.231 & 0.208 & 0.197\\
    BootLasso & 0.118 & 0.061 & 0.171 & 0.147 & 0.169 \\
  PCalg & 0.116 & 0.089 & 0.171 & 0.149 & 0.130 \\
  GSE & 0.101 & 0.089& 0.170 & 0.153 &0.133\\
\noalign{\smallskip}\hline
 \end{tabular}
 \end{center}
\end{table}

\section{Conclusion and discussion}
\label{sec:conc}
In this paper, we proposed a hybrid genetic/convex algorithm for inferring large graphs
based on a particular decomposition of the $\ell_1$-penalized maximum likelihood criterion.
We obtained two convergence inequalities that ensure that the graph estimator converges to the true graph under assumptions that mainly control the model structure: graph size
(balance between sparsity, number of nodes and maximal degree) and signal-to-noise ratio.
From an algorithmic point of view, the estimation task is split into two subproblems: 
node ordering estimation and graph structure learning. 
The first one is a non-trivial problem since we optimize over a discrete non-convex large dimensional set. It led us to use a heuristic approach we specifically tailored to achieve the optimization task.
The second one is a more common problem, related to the Lasso one, for which we proposed a sound procedure with theoretical guarantees.
The potential of such an approach clearly appeared in the numerical experiments, for which
the behavior of our algorithm seemed to be very competitive when compared to the state-of-the-art.

Nevertheless, we see many opportunities for further improvements. 
First, convergence proof for the algorithm, although a challenging task, is worth investigating, for instance using
the works of \citet{Cerf98} on genetic algorithms. An alternative would be to consider other optimization schemes for the node ordering
with more established convergence proofs (\textit{e.g.} simulated annealing \citep{granville1994}).

Second, other potential extensions involve algorithmic considerations in order to improve the calculation time, 
including a finer calibration of the algorithm parameters, an initialization step for the gradient descent, 
and, in general, improving the interactions between the nested and outer loops. Tackling very large datasets
from several thousands of nodes may also require a particular treatment, for instance by adding local search
operators to GADAG.

Finally, we would like to emphasize the graph identifiability problem: in our settings, we assume the noise variances 
of all graph nodes to be equal to ensure graph identifiability (that is no equivalence class of graphs). 
Such a hypothesis is of course restrictive and likely to be violated for real datasets.
In order to infer networks for any noise variances, one solution consists in incorporating interventional data on the model.
These data are obtained from perturbations of the biological system (\textit{e.g.} gene knockouts or over-expressions) and make the equivalence class of graphs smaller \citep{Hauser12}. 
The use of additional data, informative yet very costly interventional data could be combined with observational on the MLE estimator. It was recently proposed by \citet{Hauser15} for a BIC-score penalized MLE, or by \citet{Rau13} for learning Gaussian Bayesian networks in the case of GRN inference. A modification of our hybrid algorithm GADAG could then lead to a more accurate identification of the true graph. Lastly, the cyclic structure framework could also be considered by using a Markov equivalence characterization \citep{Richardson97} to relax the strictly triangular assumption on our matrix $T$ using Equation of Proposition \ref{prop:decomp}. It would pave the way for totally new theoretical developments, and a more realistic modelling of a complex system.


\bibliographystyle{biometrika}
\bibliography{refs-3}

\end{document}